\documentclass[11pt]{amsart}
\usepackage{amsmath,amssymb,stmaryrd}
\newtheorem{theorem}{Theorem}
\newtheorem{proposition}[theorem]{Proposition}
\newtheorem{corollary}[theorem]{Corollary}
\newtheorem{lemma}[theorem]{Lemma}

\begin{document}

\title[Harnack inequality for the Ricci flow]{A generalization of Hamilton's differential Harnack inequality for the Ricci flow}
\author{Simon Brendle}
\address{Department of Mathematics \\
                 Stanford University \\
                 Stanford, CA 94305}
\thanks{This project was supported by the Alfred P. Sloan Foundation and by the National Science Foundation under grant DMS-0605223.}
\maketitle

\section{Introduction} 

In \cite{Hamilton4}, R.~Hamilton established a differential Harnack inequality for solutions to the Ricci flow with nonnegative curvature operator (see \cite{Hamilton3} for an earlier result in the two-dimensional case). This inequality has since become one of the fundamental tools in the study of Ricci flow. We point out that H.D.~Cao \cite{Cao} has proved a differential Harnack inequality for solutions to the K\"ahler-Ricci flow with nonnegative holomorphic bisectional curvature.

In this paper, we prove a generalization of Hamilton's Harnack inequality, replacing the assumption of nonnegative curvature operator by a weaker curvature condition. Throughout this paper, we assume that $(M,g(t))$, $t \in (0,T)$, is a family of complete Riemannian manifolds evolving under Ricci flow. Following R.~Hamilton \cite{Hamilton4}, we define 
\[P_{ijk} = D_i \text{\rm Ric}_{jk} - D_j \text{\rm Ric}_{ik}\] 
and 
\[M_{ij} = \Delta \text{\rm Ric}_{ij} - \frac{1}{2} \, D_i D_j \text{\rm scal} + 2 \, R_{ikjl}  \, \text{\rm Ric}^{kl} - \text{\rm Ric}_i^k \, \text{\rm Ric}_{jk} + \frac{1}{2t} \, \text{\rm Ric}_{ij}.\] 
Here, $\text{\rm Ric}$ and $\text{\rm scal}$ denote the Ricci and scalar curvature of $(M,g(t))$, respectively. 

\begin{theorem}
\label{Harnack.inequality}
Suppose that $(M,g(t)) \times \mathbb{R}^2$ has nonnegative isotropic curvature for all $t \in (0,T)$. Moreover, we assume that 
\[\sup_{(x,t) \in M \times (\alpha,T)} \text{\rm scal}(x,t) < \infty\] 
for all $\alpha \in (0,T)$. Then 
\[M(w,w) + 2 \, P(v,w,w) + R(v,w,v,w) \geq 0\] 
for all points $(x,t) \in M \times (0,T)$ and all vectors $v,w \in T_x M$.
\end{theorem}

As a consequence, we obtain a generalization of Hamilton's trace Harnack inequality (cf. \cite{Hamilton4}):

\begin{corollary} 
\label{trace.Harnack.inequality}
Suppose that $(M,g(t)) \times \mathbb{R}^2$ has nonnegative isotropic curvature for all $t \in (0,T)$. Moreover, we assume that 
\[\sup_{(x,t) \in M \times (\alpha,T)} \text{\rm scal}(x,t) < \infty\] 
for all $\alpha \in (0,T)$. Then we have 
\[\frac{\partial}{\partial t} \text{\rm scal} + \frac{1}{t} \, \text{\rm scal} + 2 \, \partial_i \text{\rm scal} \, v^i + 2 \, \text{\rm Ric}(v,v) \geq 0\] 
for all points $(x,t) \in M \times (0,T)$ and all vectors $v \in T_x M$. 
\end{corollary}

The condition that $M \times \mathbb{R}^2$ has nonnegative isotropic curvature is preserved by the Ricci flow, and 
plays a key role in the proof of the $1/4$-pinching theorem \cite{Brendle-Schoen1}. We point out that the following statements are equivalent: 
\begin{itemize}
\item[(i)] The product $M \times \mathbb{R}^2$ has nonnegative isotropic curvature.
\item[(ii)] For all orthonormal four-frames $\{e_1,e_2,e_3,e_4\} \subset T_x M$ and all $\lambda,\mu \in [-1,1]$, we have 
\begin{align*} 
&R(e_1,e_3,e_1,e_3) + \lambda^2 \, R(e_1,e_4,e_1,e_4) \\ 
&+ \mu^2 \, R(e_2,e_3,e_2,e_3) + \lambda^2\mu^2 \, R(e_2,e_4,e_2,e_4) \\ &- 2\lambda\mu \, R(e_1,e_2,e_3,e_4) \geq 0. 
\end{align*} 
\item[(iii)] For all vectors $v_1,v_2,v_3,v_4 \in T_x M$, we have
\begin{align*} 
&R(v_1,v_3,v_1,v_3) + R(v_1,v_4,v_1,v_4) \\ 
&+ R(v_2,v_3,v_2,v_3) + R(v_2,v_4,v_2,v_4) \\ &- 2 \, R(v_1,v_2,v_3,v_4) \geq 0. 
\end{align*} 
\end{itemize}
The implication $\text{\rm (i)} \Longrightarrow \text{\rm (ii)}$ was established in \cite{Brendle-Schoen1}. Moreover, a careful examination of the proof of Proposition 21 in \cite{Brendle-Schoen1} shows that (ii) implies (iii). Finally, the implication $\text{\rm (iii)} \Longrightarrow \text{\rm (i)}$ is trivial.

\section{The space-time curvature tensor and its evolution under Ricci flow}

We first review the evolution equation for the various quantities that appear in the Harnack 
inequaltiy. The evolution equation of the curvature tensor is given by 
\begin{align*} 
&\frac{\partial}{\partial t} R_{ijkl} - \Delta R_{ijkl} \\ 
&+ \text{\rm Ric}_i^m \, R_{mjkl} + \text{\rm Ric}_j^m \, R_{imkl} + \text{\rm Ric}_k^m \, R_{ijml} + \text{\rm Ric}_l^m \, R_{ijkm} \\ 
&= g^{pq} \, g^{rs} \, R_{ijpr} \, R_{klqs} + 2 \, g^{pq} \, g^{rs} \, R_{ipkr} \, R_{iqls} - 2 \, g^{pq} \, g^{rs} \, R_{iplr} \, R_{jpks} 
\end{align*} 
(cf. \cite{Hamilton1},\cite{Hamilton2}). Moreover, Hamilton proved that 
\begin{align*} 
&\frac{\partial}{\partial t} P_{ijk} - \Delta P_{ijk} + \text{\rm Ric}_i^m \, P_{mjk} + \text{\rm Ric}_j^m \, P_{imk} + \text{\rm Ric}_k^m \, P_{ijm} \notag \\ 
&= 2 \, g^{pq} \, g^{rs} \, R_{ipjr} \, P_{qsk} + 2 \, g^{pq} \, g^{rs} \, R_{ipkr} \, P_{qjs} + 2 \, g^{pq} \, g^{rs} \, R_{jpkr} \, P_{iqs} \\ 
&- 2 \, \text{\rm Ric}_p^m \, D^p R_{ijkm} 
\end{align*} 
and 
\begin{align*} 
&\frac{\partial}{\partial t} M_{ij} - \Delta M_{ij} + \text{\rm Ric}_i^m \, M_{mj} + \text{\rm Ric}_j^m \, M_{im} \\ 
&= 2 \, g^{pq} \, g^{rs} \, R_{ipjr} \, M_{qs} + 2 \, \text{\rm Ric}_p^m \, (D^p P_{mij} + D^p P_{mji}) \\ 
&+ 2 \, g^{pq} \, g^{rs} \, P_{ipr} \, P_{jqs} - 4 \, g^{pq} \, g^{rs} \, P_{ipr} \, P_{jsq} + 2 \, \text{\rm Ric}^{lp} \, \text{\rm Ric}_p^m \, R_{iljm} - \frac{1}{2t^2} \, \text{\rm Ric}_{ij} 
\end{align*} 
(see \cite{Hamilton4}, Lemma 4.3 and Lemma 4.4). 

Chow and Chu \cite{Chow-Chu} observed that the quantities $M_{ij}$ and $P_{ijk}$ can be viewed as components of a space-time curvature tensor (see also \cite{Chow-Knopf}). In the remainder of this section, we describe the definition of the space-time curvature tensor, and its evolution under Ricci flow. Following \cite{Chow-Knopf}, we define a connection $\tilde{D}$ on the product $M \times (0,T)$ by 
\begin{align*} 
&\tilde{D}_{\frac{\partial}{\partial x^i}} \frac{\partial}{\partial x^j} = \Gamma_{ij}^k \, \frac{\partial}{\partial x^k} \\ 
&\tilde{D}_{\frac{\partial}{\partial x^i}} \frac{\partial}{\partial t} = -\Big ( \text{\rm Ric}_i^j + \frac{1}{2t} \, \delta_i^j \Big ) \, \frac{\partial}{\partial x^j} \\ 
&\tilde{D}_{\frac{\partial}{\partial t}} \frac{\partial}{\partial x^i} = -\Big ( \text{\rm Ric}_i^j + \frac{1}{2t} \, \delta_i^j \Big ) \, \frac{\partial}{\partial x^j} \\ 
&\tilde{D}_{\frac{\partial}{\partial t}} \frac{\partial}{\partial t} = -\frac{1}{2} \, \partial^i \text{\rm scal} \, \frac{\partial}{\partial x^i} - \frac{3}{2t} \, \frac{\partial}{\partial t}. 
\end{align*} 
Here, $\Gamma_{ij}^k$ denote the Christoffel symbols associated with the metric $g(t)$. We next define a $(0,4)$-tensor $S$ by 
\begin{align*} 
S &= R_{ijkl} \, dx^i \otimes dx^j \otimes dx^k \otimes dx^l \\ 
&+ P_{ijk} \, dx^i \otimes dx^j \otimes dt \otimes dx^k - P_{ijk} \, dx^i \otimes dx^j \otimes dx^k \otimes dt \\ 
&+ P_{ijk} \, dt \otimes dx^k \otimes dx^i \otimes dx^j - P_{ijk} \, dx^k \otimes dt \otimes dx^i \otimes dx^j \\ 
&+ M_{ij} \, dx^i \otimes dt \otimes dx^j \otimes dt - M_{ij} \, dx^i \otimes dt \otimes dt \otimes dx^j \\ 
&- M_{ij} \, dt \otimes dx^i \otimes dx^j \otimes dt + M_{ij} \, dt \otimes dx^i \otimes dt \otimes dx^j. 
\end{align*} 
The tensor $S$ is an algebraic curvature tensor in the sense that 
\[S(\tilde{v}_1,\tilde{v}_2,\tilde{v}_3,\tilde{v}_4) = -S(\tilde{v}_2,\tilde{v}_1,\tilde{v}_3,\tilde{v}_4) = S(\tilde{v}_3,\tilde{v}_4,\tilde{v}_1,\tilde{v}_2)\] 
and 
\[S(\tilde{v}_1,\tilde{v}_2,\tilde{v}_3,\tilde{v}_4) + 
S(\tilde{v}_2,\tilde{v}_3,\tilde{v}_1,\tilde{v}_4) + S(\tilde{v}_3,\tilde{v}_1,\tilde{v}_2,\tilde{v}_4) = 0\] 
for all vectors $\tilde{v}_1,\tilde{v}_2,\tilde{v}_3,\tilde{v}_4 \in T_{(x,t)}(M \times (0,T))$.

Given any algebraic curvature tensor $S$, we define 
\begin{align*} 
&\tilde{Q}(S)(\tilde{v}_1,\tilde{v}_2,\tilde{v}_3,\tilde{v}_4) \\ 
&= \sum_{p,q,r,s=1}^n g^{pq} \, g^{rs} \, S \Big ( \tilde{v}_1,\tilde{v}_2,\frac{\partial}{\partial x^p},\frac{\partial}{\partial x^r} \Big ) \, S \Big ( \tilde{v}_3,\tilde{v}_4,\frac{\partial}{\partial x^q},\frac{\partial}{\partial x^s} \Big ) \\ 
&+ 2 \sum_{p,q,r,s=1}^n g^{pq} \, g^{rs} \, S \Big ( \tilde{v}_1,\frac{\partial}{\partial x^p},\tilde{v}_3,\frac{\partial}{\partial x^r} \Big ) \, S \Big ( \tilde{v}_2,\frac{\partial}{\partial x^q},\tilde{v}_4,\frac{\partial}{\partial x^s} \Big ) \\ 
&- 2 \sum_{p,q,r,s=1}^n g^{pq} \, g^{rs} \, S \Big ( \tilde{v}_1,\frac{\partial}{\partial x^p},\tilde{v}_4,\frac{\partial}{\partial x^r} \Big ) \, S \Big ( \tilde{v}_2,\frac{\partial}{\partial x^q},\tilde{v}_3,\frac{\partial}{\partial x^s} \Big ) 
\end{align*} 
for all vectors $\tilde{v}_1,\tilde{v}_2,\tilde{v}_3,\tilde{v}_4 \in T_{(x,t)} (M \times (0,T))$. 
It is straightforward to verify that 
\[\tilde{Q}(S)(\tilde{v}_1,\tilde{v}_2,\tilde{v}_3,\tilde{v}_4) = -\tilde{Q}(S)(\tilde{v}_2,\tilde{v}_1,\tilde{v}_3,\tilde{v}_4) = \tilde{Q}(S)(\tilde{v}_3,\tilde{v}_4,\tilde{v}_1,\tilde{v}_2)\] 
and 
\[\tilde{Q}(S)(\tilde{v}_1,\tilde{v}_2,\tilde{v}_3,\tilde{v}_4) + 
\tilde{Q}(S)(\tilde{v}_2,\tilde{v}_3,\tilde{v}_1,\tilde{v}_4) + \tilde{Q}(S)(\tilde{v}_3,\tilde{v}_1,\tilde{v}_2,\tilde{v}_4) = 0\] 
for all vectors $\tilde{v}_1,\tilde{v}_2,\tilde{v}_3,\tilde{v}_4 \in T_{(x,t)}(M \times (0,T))$.
Therefore, $\tilde{Q}(S)$ is again an algebraic curvature tensor. \\

\begin{proposition} 
\label{evolution.equation}
The tensor $S$ satisfies the evolution equation 
\[\tilde{D}_{\frac{\partial}{\partial t}} S = \tilde{\Delta} S + \frac{2}{t} \, S + \tilde{Q}(S).\] 
Here, 
\[\tilde{\Delta} S = \sum_{p,q=1}^n g^{pq} \, \tilde{D}_{\frac{\partial}{\partial x^p},\frac{\partial}{\partial x^q}}^2 S\] 
denotes the Laplacian of $S$ with respect to the connection $\tilde{D}$.
\end{proposition}

\textbf{Proof.} 
For abbreviation, let $W = \tilde{D}_{\frac{\partial}{\partial t}} S - \tilde{\Delta} S - \frac{2}{t} \, S$. Clearly, $W$ is an algebraic curvature tensor. 
We claim that $W = \tilde{Q}(S)$. Note that 
\begin{align*} 
&(\tilde{D}_{\frac{\partial}{\partial t}} S) \Big ( \frac{\partial}{\partial x^i},\frac{\partial}{\partial x^j},\frac{\partial}{\partial x^k},\frac{\partial}{\partial x^l} \Big ) \\ 
&= \frac{\partial}{\partial t} R_{ijkl} + \frac{2}{t} \, R_{ijkl} \\ &+ \text{\rm Ric}_i^m \, R_{mjkl} + \text{\rm Ric}_j^m \, R_{imkl} + \text{\rm Ric}_k^m \, R_{ijml} + \text{\rm Ric}_l^m \, R_{ijkm}  
\end{align*} 
and 
\[(\tilde{\Delta} S) \Big ( \frac{\partial}{\partial x^i},\frac{\partial}{\partial x^j},\frac{\partial}{\partial x^k},\frac{\partial}{\partial x^l} \Big ) = \Delta R_{ijkl}.\] 
This implies 
\begin{align*} 
&W \Big ( \frac{\partial}{\partial x^i},\frac{\partial}{\partial x^j},\frac{\partial}{\partial x^k},\frac{\partial}{\partial x^l} \Big ) \\ 
&= \frac{\partial}{\partial t} R_{ijkl} - \Delta R_{ijkl} \\ &+ \text{\rm Ric}_i^m \, R_{mjkl} + \text{\rm Ric}_j^m \, R_{imkl} + \text{\rm Ric}_k^m \, R_{ijml} + \text{\rm Ric}_l^m \, R_{ijkm} \\ 
&= g^{pq} \, g^{rs} \, R_{ijpr} \, R_{klqs} + 2 \, g^{pq} \, g^{rs} \, R_{ipkr} \, R_{iqls} - 2 \, g^{pq} \, g^{rs} \, R_{iplr} \, R_{jpks} \\ 
&= \tilde{Q}(S) \Big ( \frac{\partial}{\partial x^i},\frac{\partial}{\partial x^j},\frac{\partial}{\partial x^k},\frac{\partial}{\partial x^l} \Big ). 
\end{align*}
Moreover, we have 
\begin{align*}
&(\tilde{D}_{\frac{\partial}{\partial t}} S) \Big ( \frac{\partial}{\partial x^i},\frac{\partial}{\partial x^j},\frac{\partial}{\partial t},\frac{\partial}{\partial x^k} \Big ) \\ 
&= \frac{\partial}{\partial t} P_{ijk} + \frac{1}{2} \, \partial^m \text{\rm scal} \, R_{ijmk} + \frac{3}{t} \, P_{ijk} \\ 
&+ \text{\rm Ric}_i^m \, P_{mjk} + \text{\rm Ric}_j^m \, P_{imk} + \text{\rm Ric}_k^m \, P_{ijm} 
\end{align*} 
and 
\begin{align*} 
&(\tilde{\Delta} S) \Big ( \frac{\partial}{\partial x^i},\frac{\partial}{\partial x^j},\frac{\partial}{\partial t},\frac{\partial}{\partial x^k} \Big ) \\ 
&= \Delta P_{ijk} + D^p \text{\rm Ric}_p^m \, R_{ijmk} + 2 \, \text{\rm Ric}_p^m \, D^p R_{ijmk} + \frac{1}{t} \, D^m R_{ijmk} \\ 
&= \Delta P_{ijk} + \frac{1}{2} \, \partial^m \text{\rm scal} \, R_{ijmk} + 2 \, \text{\rm Ric}_p^m \, D^p R_{ijmk} + \frac{1}{t} \, P_{ijk}. 
\end{align*}
Using the evolution equation for the tensor $P_{ijk}$, we obtain 
\begin{align*}
&W \Big ( \frac{\partial}{\partial x^i},\frac{\partial}{\partial x^j},\frac{\partial}{\partial t},\frac{\partial}{\partial x^k} \Big ) \\
&= \frac{\partial}{\partial t} P_{ijk} - \Delta P_{ijk} - 2 \, \text{\rm Ric}_p^m \, D^p R_{ijmk} \\ &+ \text{\rm Ric}_i^m \, P_{mjk} + \text{\rm Ric}_j^m \, P_{imk} + \text{\rm Ric}_k^m \, P_{ijm} \\ 
&= 2 \, g^{pq} \, g^{rs} \, R_{ipjr} \, P_{qsk} + 2 \, g^{pq} \, g^{rs} \, R_{ipkr} \, P_{qjs} + 2 \, g^{pq} \, g^{rs} \, R_{jpkr} \, P_{iqs} \\ 
&= g^{pq} \, g^{rs} \, R_{ijpr} \, P_{qsk} + 2 \, g^{pq} \, g^{rs} \, R_{ipkr} \, P_{qjs} + 2 \, g^{pq} \, g^{rs} \, R_{jpkr} \, P_{iqs} \\ 
&= \tilde{Q}(S) \Big ( \frac{\partial}{\partial x^i},\frac{\partial}{\partial x^j},\frac{\partial}{\partial t},\frac{\partial}{\partial x^k} \Big ). 
\end{align*} 
Finally, we have 
\begin{align*}
&(\tilde{D}_{\frac{\partial}{\partial t}} S) \Big ( \frac{\partial}{\partial x^i},\frac{\partial}{\partial t},\frac{\partial}{\partial x^j},\frac{\partial}{\partial t} \Big ) \\ 
&= \frac{\partial}{\partial t} M_{ij} - \frac{1}{2} \, \partial^m \text{\rm scal} \, (P_{imj} + P_{jmi}) + \frac{4}{t} \, M_{ij} + \text{\rm Ric}_i^m \, M_{mj} + \text{\rm Ric}_j^m \, M_{im} 
\end{align*} 
and 
\begin{align*} 
&(\tilde{\Delta} S) \Big ( \frac{\partial}{\partial x^i},\frac{\partial}{\partial t},\frac{\partial}{\partial x^j},\frac{\partial}{\partial t} \Big ) \\ 
&= \Delta M_{ij} - 2 \, \Big ( \text{\rm Ric}_p^m + \frac{1}{2t} \, \delta_p^m \Big ) \, (D^p P_{imj} + D^p P_{jmi}) \\ 
&- D^p \text{\rm Ric}_p^m \, (P_{imj} + P_{jmi}) + 2 \, \Big ( \text{\rm Ric}^{lp} + \frac{1}{2t} \, g^{lp} \Big ) \, 
\Big ( \text{\rm Ric}_p^m + \frac{1}{2t} \, \delta_p^m \Big ) \, R_{iljm} \\ 
&= \Delta M_{ij} - 2 \, \text{\rm Ric}_p^m \, (D^p P_{imj} + D^p P_{jmi}) \\ 
&- D^p \text{\rm Ric}_p^m \, (P_{imj} + P_{jmi}) + 2 \, \text{\rm Ric}^{lp} \, \text{\rm Ric}_p^m \, R_{iljm} \\ 
&- \frac{1}{t} \, (D^m P_{imj} + D^m P_{jmi}) + \frac{2}{t} \, \text{\rm Ric}^{lm} \, R_{iljm} 
+ \frac{1}{2t^2} \, \text{\rm Ric}_{ij} \\
&= \Delta M_{ij} - 2 \, \text{\rm Ric}_p^m \, (D^p P_{imj} + D^p P_{jmi}) \\ 
&- \frac{1}{2} \, \partial^m \text{\rm scal} \, (P_{imj} + P_{jmi}) + 2 \, \text{\rm Ric}^{lp} \, \text{\rm Ric}_p^m \, R_{iljm} \\ 
&+ \frac{2}{t} \, M_{ij} - \frac{1}{2t^2} \, \text{\rm Ric}_{ij}. 
\end{align*} 
In the last step, we have used the formula 
\[M_{ij} = -D^m P_{imj} + \text{\rm Ric}^{lm} \, R_{iljm} + \frac{1}{2t} \, \text{\rm Ric}_{ij}\] 
(see \cite{Hamilton4}, p.~235). Using the evolution equation for $M_{ij}$, we obtain 
\begin{align*}
&W \Big ( \frac{\partial}{\partial x^i},\frac{\partial}{\partial t},\frac{\partial}{\partial x^j},\frac{\partial}{\partial t} \Big ) \\ 
&= \frac{\partial}{\partial t} M_{ij} - \Delta M_{ij} + 2 \, \text{\rm Ric}_p^m \, (D^p P_{imj} + D^p P_{jmi}) \\ 
&+ \text{\rm Ric}_i^m \, M_{mj} + \text{\rm Ric}_j^m \, M_{im} - 2 \, \text{\rm Ric}^{lp} \, \text{\rm Ric}_p^m \, R_{iljm} + \frac{1}{2t^2} \, \text{\rm Ric}_{ij} \\ 
&= 2 \, g^{pq} \, g^{rs} \, R_{ipjr} \, M_{qs} + 2 \, g^{pq} \, g^{rs} \, P_{ipr} \, P_{jqs} - 4 \, g^{pq} \, g^{rs} \, P_{ipr} \, P_{jsq} \\ 
&= 2 \, g^{pq} \, g^{rs} \, R_{ipjr} \, M_{qs} + 2 \, g^{pq} \, g^{rs} \, P_{ipr} \, (P_{jqs} - P_{jsq}) - 2 \, g^{pq} \, g^{rs} \, P_{ipr} \, P_{jsq} \\ 
&= 2 \, g^{pq} \, g^{rs} \, R_{ipjr} \, M_{qs} - 2 \, g^{pq} \, g^{rs} \, P_{ipr} \, P_{qsj} - 2 \, g^{pq} \, g^{rs} \, P_{ipr} \, P_{jsq} \\ 
&= 2 \, g^{pq} \, g^{rs} \, R_{ipjr} \, M_{qs} + g^{pq} \, g^{rs} \, P_{pri} \, P_{qsj} - 2 \, g^{pq} \, g^{rs} \, P_{ipr} \, P_{jsq} \\ 
&= \tilde{Q}(S) \Big ( \frac{\partial}{\partial x^i},\frac{\partial}{\partial t},\frac{\partial}{\partial x^j},\frac{\partial}{\partial t} \Big ). 
\end{align*} 
Putting these facts together, we conclude that $W = \tilde{Q}(S)$. This completes the proof.

\section{An invariant cone for the ODE $\frac{d}{dt} S = \tilde{Q}(S)$}

We now consider the space of algebraic curvature tensors on $\mathbb{R}^n \times \mathbb{R}$. There is a natural mapping $\tilde{Q}$ which maps the space of algebraic curvature tensors on $\mathbb{R}^n \times \mathbb{R}$ into itself. For each algebraic curvature tensor $S$ on $\mathbb{R}^n \times \mathbb{R}$, the tensor $\tilde{Q}(S)$ is defined by 
\begin{align*} 
\tilde{Q}(S)(\tilde{v}_1,\tilde{v}_2,\tilde{v}_3,\tilde{v}_4) 
&= \sum_{p,q=1}^n S(\tilde{v}_1,\tilde{v}_2,e_p,e_q) \, S(\tilde{v}_3,\tilde{v}_4,e_p,e_q) \\ 
&+ 2 \sum_{p,q=1}^n S(\tilde{v}_1,e_p,\tilde{v}_3,e_q) \, S(\tilde{v}_2,e_p,\tilde{v}_4,e_q) \\ 
&- 2 \sum_{p,q=1}^n S(\tilde{v}_1,e_p,\tilde{v}_4,e_q) \, S(\tilde{v}_2,e_p,\tilde{v}_3,e_q), 
\end{align*} 
where $\{e_1,\hdots,e_n\}$ is an orthonormal basis of $\mathbb{R}^n$. 

Let $K$ be the set of all algebraic curvature tensors on $\mathbb{R}^n \times \mathbb{R}$ such that 
\begin{align*} 
&S(\tilde{v}_1,\tilde{v}_3,\tilde{v}_1,\tilde{v}_3) + S(\tilde{v}_1,\tilde{v}_4,\tilde{v}_1,\tilde{v}_4) \\ 
&+ S(\tilde{v}_2,\tilde{v}_3,\tilde{v}_2,\tilde{v}_3) + S(\tilde{v}_2,\tilde{v}_4,\tilde{v}_2,\tilde{v}_4) - 2 \, S(\tilde{v}_1,\tilde{v}_2,\tilde{v}_3,\tilde{v}_4) \geq 0 
\end{align*} 
for all vectors $\tilde{v}_1,\tilde{v}_2,\tilde{v}_3,\tilde{v}_4 \in \mathbb{R}^n \times \mathbb{R}$. Clearly, $K$ is a closed convex cone. Moreover, $K$ is invariant under the natural action of $GL(n+1)$. 

We claim that $K$ is invariant under the ODE $\frac{d}{dt} S = \tilde{Q}(S)$. The proof relies on the following result:

\begin{proposition}
\label{second.variation}
Let $S$ be an algebraic curvature tensor on $\mathbb{R}^n \times \mathbb{R}$ which lies in the cone $K$. Moreover, suppose that $\tilde{v}_1,\tilde{v}_2,\tilde{v}_3,\tilde{v}_4$ are vectors in $\mathbb{R}^n \times \mathbb{R}$ satisfying 
\begin{align*} 
&S(\tilde{v}_1,\tilde{v}_3,\tilde{v}_1,\tilde{v}_3) + S(\tilde{v}_1,\tilde{v}_4,\tilde{v}_1,\tilde{v}_4) \\ 
&+ S(\tilde{v}_2,\tilde{v}_3,\tilde{v}_2,\tilde{v}_3) + S(\tilde{v}_2,\tilde{v}_4,\tilde{v}_2,\tilde{v}_4) - 2 \, S(\tilde{v}_1,\tilde{v}_2,\tilde{v}_3,\tilde{v}_4) = 0. 
\end{align*} 
Then the expression 
\begin{align*}
&S(\tilde{w}_1,\tilde{v}_3,\tilde{w}_1,\tilde{v}_3) + S(\tilde{w}_1,\tilde{v}_4,\tilde{w}_1,\tilde{v}_4) \\ 
&+ S(\tilde{w}_2,\tilde{v}_3,\tilde{w}_2,\tilde{v}_3) + S(\tilde{w}_2,\tilde{v}_4,\tilde{w}_2,\tilde{v}_4) \\ 
&+ S(\tilde{v}_1,\tilde{w}_3,\tilde{v}_1,\tilde{w}_3) + S(\tilde{v}_2,\tilde{w}_3,\tilde{v}_2,\tilde{w}_3) \\ 
&+ S(\tilde{v}_1,\tilde{w}_4,\tilde{v}_1,\tilde{w}_4) + S(\tilde{v}_2,\tilde{w}_4,\tilde{v}_2,\tilde{w}_4) \\ 
&- 2 \, \big [ S(\tilde{v}_3,\tilde{w}_1,\tilde{v}_1,\tilde{w}_3) + S(\tilde{v}_4,\tilde{w}_1,\tilde{v}_2,\tilde{w}_3) \big ] \\ 
&- 2 \, \big [ S(\tilde{v}_4,\tilde{w}_1,\tilde{v}_1,\tilde{w}_4) - S(\tilde{v}_3,\tilde{w}_1,\tilde{v}_2,\tilde{w}_4) \big ] \\ 
&+ 2 \, \big [ S(\tilde{v}_4,\tilde{w}_2,\tilde{v}_1,\tilde{w}_3) - S(\tilde{v}_3,\tilde{w}_2,\tilde{v}_2,\tilde{w}_3) \big ] \\ 
&- 2 \, \big [ S(\tilde{v}_3,\tilde{w}_2,\tilde{v}_1,\tilde{w}_4) + S(\tilde{v}_4,\tilde{w}_2,\tilde{v}_2,\tilde{w}_4) \big ] \\ 
&- 2 \, S(\tilde{w}_1,\tilde{w}_2,\tilde{v}_3,\tilde{v}_4) - 2 \, S(\tilde{v}_1,\tilde{v}_2,\tilde{w}_3,\tilde{w}_4) 
\end{align*} 
is nonnegative for all vectors $\tilde{w}_1,\tilde{w}_2,\tilde{w}_3,\tilde{w}_4 \in \mathbb{R}^n \times \mathbb{R}$.
\end{proposition}

\textbf{Proof.} 
The proof is similar to the proof of Proposition 8 in \cite{Brendle-Schoen1}. Since $S \in K$, we have 
\begin{align*} 
0 &\leq S(\tilde{v}_1+s\tilde{w}_1,\tilde{v}_3+s\tilde{w}_3,\tilde{v}_1+s\tilde{w}_1,\tilde{v}_3+s\tilde{w}_3) \\ 
&+ S(\tilde{v}_1+s\tilde{w}_1,\tilde{v}_4+s\tilde{w}_4,\tilde{v}_1+s\tilde{w}_1,\tilde{v}_4+s\tilde{w}_4) \\ 
&+ S(\tilde{v}_2+s\tilde{w}_2,\tilde{v}_3+s\tilde{w}_3,\tilde{v}_2+s\tilde{w}_2,\tilde{v}_3+s\tilde{w}_3) \\ 
&+ S(\tilde{v}_2+s\tilde{w}_2,\tilde{v}_4+s\tilde{w}_4,\tilde{v}_2+s\tilde{w}_2,\tilde{v}_4+s\tilde{w}_4) \\ 
&- 2 \, S(\tilde{v}_1+s\tilde{w}_1,\tilde{v}_2+s\tilde{w}_2,\tilde{v}_3+s\tilde{w}_3,\tilde{v}_4+s\tilde{w}_4) 
\end{align*} 
for all $s \in \mathbb{R}$. Taking the second derivative at $s = 0$, we obtain 
\begin{align}
\label{inequality.1}
0 &\leq S(\tilde{w}_1,\tilde{v}_3,\tilde{w}_1,\tilde{v}_3) + S(\tilde{w}_1,\tilde{v}_4,\tilde{w}_1,\tilde{v}_4) \notag \\ 
&+ S(\tilde{w}_2,\tilde{v}_3,\tilde{w}_2,\tilde{v}_3) + S(\tilde{w}_2,\tilde{v}_4,\tilde{w}_2,\tilde{v}_4) \notag \\ 
&+ S(\tilde{v}_1,\tilde{w}_3,\tilde{v}_1,\tilde{w}_3) + S(\tilde{v}_2,\tilde{w}_3,\tilde{v}_2,\tilde{w}_3) \notag \\ 
&+ S(\tilde{v}_1,\tilde{w}_4,\tilde{v}_1,\tilde{w}_4) + S(\tilde{v}_2,\tilde{w}_4,\tilde{v}_2,\tilde{w}_4) \notag \\ 
&+ 2 \, S(\tilde{v}_1,\tilde{v}_3,\tilde{w}_1,\tilde{w}_3) + 2 \, S(\tilde{v}_1,\tilde{w}_3,\tilde{w}_1,\tilde{v}_3) - 2 \, S(\tilde{w}_1,\tilde{v}_2,\tilde{w}_3,\tilde{v}_4) \\ 
&+ 2 \, S(\tilde{v}_1,\tilde{v}_4,\tilde{w}_1,\tilde{w}_4) + 2 \, S(\tilde{v}_1,\tilde{w}_4,\tilde{w}_1,\tilde{v}_4) - 2 \, S(\tilde{w}_1,\tilde{v}_2,\tilde{v}_3,\tilde{w}_4) \notag \\ 
&+ 2 \, S(\tilde{v}_2,\tilde{v}_3,\tilde{w}_2,\tilde{w}_3) + 2 \, S(\tilde{v}_2,\tilde{w}_3,\tilde{w}_2,\tilde{v}_3) - 2 \, S(\tilde{v}_1,\tilde{w}_2,\tilde{w}_3,\tilde{v}_4) \notag \\ 
&+ 2 \, S(\tilde{v}_2,\tilde{v}_4,\tilde{w}_2,\tilde{w}_4) + 2 \, S(\tilde{v}_2,\tilde{w}_4,\tilde{w}_2,\tilde{v}_4) - 2 \, S(\tilde{v}_1,\tilde{w}_2,\tilde{v}_3,\tilde{w}_4) \notag \\ 
&- 2 \, S(\tilde{w}_1,\tilde{w}_2,\tilde{v}_3,\tilde{v}_4) - 2 \, S(\tilde{v}_1,\tilde{v}_2,\tilde{w}_3,\tilde{w}_4). \notag 
\end{align}
Replacing $\{\tilde{v}_1,\tilde{v}_2,\tilde{v}_3,\tilde{v}_4\}$ by $\{\tilde{v}_2,-\tilde{v}_1,\tilde{v}_4,-\tilde{v}_3\}$ yields 
\begin{align}
\label{inequality.2}
0 &\leq S(\tilde{w}_1,\tilde{v}_4,\tilde{w}_1,\tilde{v}_4) + S(\tilde{w}_1,\tilde{v}_3,\tilde{w}_1,\tilde{v}_3) \notag \\ 
&+ S(\tilde{w}_2,\tilde{v}_4,\tilde{w}_2,\tilde{v}_4) + S(\tilde{w}_2,\tilde{v}_3,\tilde{w}_2,\tilde{v}_3) \notag \\ 
&+ S(\tilde{v}_2,\tilde{w}_3,\tilde{v}_2,\tilde{w}_3) + S(\tilde{v}_1,\tilde{w}_3,\tilde{v}_1,\tilde{w}_3) \notag \\ 
&+ S(\tilde{v}_2,\tilde{w}_4,\tilde{v}_2,\tilde{w}_4) + S(\tilde{v}_1,\tilde{w}_4,\tilde{v}_1,\tilde{w}_4) \notag \\ 
&+ 2 \, S(\tilde{v}_2,\tilde{v}_4,\tilde{w}_1,\tilde{w}_3) + 2 \, S(\tilde{v}_2,\tilde{w}_3,\tilde{w}_1,\tilde{v}_4) - 2 \, S(\tilde{w}_1,\tilde{v}_1,\tilde{w}_3,\tilde{v}_3) \\ 
&- 2 \, S(\tilde{v}_2,\tilde{v}_3,\tilde{w}_1,\tilde{w}_4) - 2 \, S(\tilde{v}_2,\tilde{w}_4,\tilde{w}_1,\tilde{v}_3) + 2 \, S(\tilde{w}_1,\tilde{v}_1,\tilde{v}_4,\tilde{w}_4) \notag \\ 
&- 2 \, S(\tilde{v}_1,\tilde{v}_4,\tilde{w}_2,\tilde{w}_3) - 2 \, S(\tilde{v}_1,\tilde{w}_3,\tilde{w}_2,\tilde{v}_4) + 2 \, S(\tilde{v}_2,\tilde{w}_2,\tilde{w}_3,\tilde{v}_3) \notag \\ 
&+ 2 \, S(\tilde{v}_1,\tilde{v}_3,\tilde{w}_2,\tilde{w}_4) + 2 \, S(\tilde{v}_1,\tilde{w}_4,\tilde{w}_2,\tilde{v}_3) - 2 \, S(\tilde{v}_2,\tilde{w}_2,\tilde{v}_4,\tilde{w}_4) \notag \\ 
&+ 2 \, S(\tilde{w}_1,\tilde{w}_2,\tilde{v}_4,\tilde{v}_3) + 2 \, S(\tilde{v}_2,\tilde{v}_1,\tilde{w}_3,\tilde{w}_4). \notag
\end{align}
In the next step, we take the arithmetic mean of (\ref{inequality.1}) and (\ref{inequality.2}). This yields 
\begin{align} 
\label{inequality.3}
0 &\leq S(\tilde{w}_1,\tilde{v}_3,\tilde{w}_1,\tilde{v}_3) + S(\tilde{w}_1,\tilde{v}_4,\tilde{w}_1,\tilde{v}_4) \notag \\ 
&+ S(\tilde{w}_2,\tilde{v}_3,\tilde{w}_2,\tilde{v}_3) + S(\tilde{w}_2,\tilde{v}_4,\tilde{w}_2,\tilde{v}_4) \notag \\ 
&+ S(\tilde{v}_1,\tilde{w}_3,\tilde{v}_1,\tilde{w}_3) + S(\tilde{v}_2,\tilde{w}_3,\tilde{v}_2,\tilde{w}_3) \notag \\ 
&+ S(\tilde{v}_1,\tilde{w}_4,\tilde{v}_1,\tilde{w}_4) + S(\tilde{v}_2,\tilde{w}_4,\tilde{v}_2,\tilde{w}_4) \notag \\ 
&+ \big [ S(\tilde{v}_1,\tilde{v}_3,\tilde{w}_1,\tilde{w}_3) + S(\tilde{v}_1,\tilde{w}_3,\tilde{w}_1,\tilde{v}_3) - S(\tilde{w}_1,\tilde{v}_2,\tilde{w}_3,\tilde{v}_4) \notag \\ 
&\hspace{2mm} + S(\tilde{v}_2,\tilde{v}_4,\tilde{w}_1,\tilde{w}_3) + S(\tilde{v}_2,\tilde{w}_3,\tilde{w}_1,\tilde{v}_4) - S(\tilde{w}_1,\tilde{v}_1,\tilde{w}_3,\tilde{v}_3) \big ] \notag \\ 
&+ \big [ S(\tilde{v}_1,\tilde{v}_4,\tilde{w}_1,\tilde{w}_4) + S(\tilde{v}_1,\tilde{w}_4,\tilde{w}_1,\tilde{v}_4) - S(\tilde{w}_1,\tilde{v}_2,\tilde{v}_3,\tilde{w}_4) \\ 
&\hspace{2mm} - S(\tilde{v}_2,\tilde{v}_3,\tilde{w}_1,\tilde{w}_4) - S(\tilde{v}_2,\tilde{w}_4,\tilde{w}_1,\tilde{v}_3) + S(\tilde{w}_1,\tilde{v}_1,\tilde{v}_4,\tilde{w}_4) \big ] \notag \\ 
&+ \big [ S(\tilde{v}_2,\tilde{v}_3,\tilde{w}_2,\tilde{w}_3) + S(\tilde{v}_2,\tilde{w}_3,\tilde{w}_2,\tilde{v}_3) - S(\tilde{v}_1,\tilde{w}_2,\tilde{w}_3,\tilde{v}_4) \notag \\ 
&\hspace{2mm} - S(\tilde{v}_1,\tilde{v}_4,\tilde{w}_2,\tilde{w}_3) - S(\tilde{v}_1,\tilde{w}_3,\tilde{w}_2,\tilde{v}_4) + S(\tilde{v}_2,\tilde{w}_2,\tilde{w}_3,\tilde{v}_3) \big ] \notag \\ 
&+ \big [ S(\tilde{v}_2,\tilde{v}_4,\tilde{w}_2,\tilde{w}_4) + S(\tilde{v}_2,\tilde{w}_4,\tilde{w}_2,\tilde{v}_4) - S(\tilde{v}_1,\tilde{w}_2,\tilde{v}_3,\tilde{w}_4) \notag \\ 
&\hspace{2mm} + S(\tilde{v}_1,\tilde{v}_3,\tilde{w}_2,\tilde{w}_4) + S(\tilde{v}_1,\tilde{w}_4,\tilde{w}_2,\tilde{v}_3) - S(\tilde{v}_2,\tilde{w}_2,\tilde{v}_4,\tilde{w}_4) \big ] \notag \\ 
&- 2 \, S(\tilde{w}_1,\tilde{w}_2,\tilde{v}_3,\tilde{v}_4) - 2 \, S(\tilde{v}_1,\tilde{v}_2,\tilde{w}_3,\tilde{w}_4). \notag 
\end{align} 
Since $S$ satisfies the first Bianchi identity, the assertion follows. \\

\begin{proposition}
\label{algebraic.fact}
Let $S$ be an algebraic curvature tensor on $\mathbb{R}^n \times \mathbb{R}$ which lies in the cone $K$. Moreover, suppose that $\tilde{v}_1,\tilde{v}_2,\tilde{v}_3,\tilde{v}_4$ are vectors in $\mathbb{R}^n \times \mathbb{R}$ satisfying 
\begin{align*} 
&S(\tilde{v}_1,\tilde{v}_3,\tilde{v}_1,\tilde{v}_3) + S(\tilde{v}_1,\tilde{v}_4,\tilde{v}_1,\tilde{v}_4) \\ 
&+ S(\tilde{v}_2,\tilde{v}_3,\tilde{v}_2,\tilde{v}_3) + S(\tilde{v}_2,\tilde{v}_4,\tilde{v}_2,\tilde{v}_4) - 2 \, S(\tilde{v}_1,\tilde{v}_2,\tilde{v}_3,\tilde{v}_4) = 0. 
\end{align*} 
Then 
\begin{align*} 
&\tilde{Q}(S)(\tilde{v}_1,\tilde{v}_3,\tilde{v}_1,\tilde{v}_3) + \tilde{Q}(S)(\tilde{v}_1,\tilde{v}_4,\tilde{v}_1,\tilde{v}_4) \\ 
&+ \tilde{Q}(S)(\tilde{v}_2,\tilde{v}_3,\tilde{v}_2,\tilde{v}_3) + \tilde{Q}(S)(\tilde{v}_2,\tilde{v}_4,\tilde{v}_2,\tilde{v}_4) 
- 2 \, \tilde{Q}(S)(\tilde{v}_1,\tilde{v}_2,\tilde{v}_3,\tilde{v}_4) \geq 0. 
\end{align*} 
\end{proposition}

\textbf{Proof.} 
Consider the following $n \times n$ matrices: 
\begin{align*} 
a_{pq} &= S(\tilde{v}_1,e_p,\tilde{v}_1,e_q) + S(\tilde{v}_2,e_p,\tilde{v}_2,e_q), \\ 
b_{pq} &= S(\tilde{v}_3,e_p,\tilde{v}_3,e_q) + S(\tilde{v}_4,e_p,\tilde{v}_4,e_q), \\ 
c_{pq} &= S(\tilde{v}_3,e_p,\tilde{v}_1,e_q) + S(\tilde{v}_4,e_p,\tilde{v}_2,e_q), \\ 
d_{pq} &= S(\tilde{v}_4,e_p,\tilde{v}_1,e_q) - S(\tilde{v}_3,e_p,\tilde{v}_2,e_q), \\ 
e_{pq} &= S(\tilde{v}_1,\tilde{v}_2,e_p,e_q), \\ 
f_{pq} &= S(\tilde{v}_3,\tilde{v}_4,e_p,e_q) 
\end{align*} 
($1 \leq p,q \leq n$). It follows from Proposition \ref{second.variation} that the matrix 
\[\begin{bmatrix} B & -F & -C & -D \\ F & B & D & -C \\ -C^T & D^T & A & -E \\ -D^T & -C^T & E & A \end{bmatrix}\] 
is positive semi-definite. This implies 
\[\text{\rm tr}(AB) + \text{\rm tr}(EF) - \text{\rm tr}(C^2) - \text{\rm tr}(D^2) \geq 0,\] 
hence 
\begin{equation} 
\label{ingredient.1}
\sum_{p,q=1}^n a_{pq} \, b_{pq} - \sum_{p,q=1}^n e_{pq} \, f_{pq} - \sum_{p,q=1}^n c_{pq} \, c_{qp} - \sum_{p,q=1}^n d_{pq} \, d_{qp} \geq 0 
\end{equation}
(cf. \cite{Brendle-Schoen1}, Proposition 9). On the other hand, we have 
\begin{align*} 
\tilde{Q}(S)(\tilde{v}_1,\tilde{v}_2,\tilde{v}_3,\tilde{v}_4) 
&= \sum_{p,q=1}^n S(\tilde{v}_1,\tilde{v}_2,e_p,e_q) \, S(\tilde{v}_3,\tilde{v}_4,e_p,e_q) \\ 
&+ \sum_{p,q=1}^n S(\tilde{v}_1,\tilde{v}_3,e_p,e_q) \, S(\tilde{v}_2,\tilde{v}_4,e_p,e_q) \\ 
&- \sum_{p,q=1}^n S(\tilde{v}_1,\tilde{v}_4,e_p,e_q) \, S(\tilde{v}_2,\tilde{v}_3,e_p,e_q) \\ 
&+ 2 \sum_{p,q=1}^n S(\tilde{v}_1,e_p,\tilde{v}_3,e_q) \, S(\tilde{v}_4,e_p,\tilde{v}_2,e_q) \\ 
&- 2 \sum_{p,q=1}^n S(\tilde{v}_1,e_p,\tilde{v}_4,e_q) \, S(\tilde{v}_3,e_p,\tilde{v}_2,e_q) 
\end{align*}
since $S$ satisfies the first Bianchi identity. This implies 
\begin{align} 
\label{ingredient.2}
&\tilde{Q}(S)(\tilde{v}_1,\tilde{v}_3,\tilde{v}_1,\tilde{v}_3) + \tilde{Q}(S)(\tilde{v}_1,\tilde{v}_4,\tilde{v}_1,\tilde{v}_4) \notag \\ 
&+ \tilde{Q}(S)(\tilde{v}_2,\tilde{v}_3,\tilde{v}_2,\tilde{v}_3) + \tilde{Q}(S)(\tilde{v}_2,\tilde{v}_4,\tilde{v}_2,\tilde{v}_4) - 2 \, \tilde{Q}(S)(\tilde{v}_1,\tilde{v}_2,\tilde{v}_3,\tilde{v}_4) \notag \\ 
&= \sum_{p,q=1}^n [S(\tilde{v}_1,\tilde{v}_3,e_p,e_q) - S(\tilde{v}_2,\tilde{v}_4,e_p,e_q)]^2 \\ 
&+ \sum_{p,q=1}^n [S(\tilde{v}_1,\tilde{v}_4,e_p,e_q) + S(\tilde{v}_2,\tilde{v}_3,e_p,e_q)]^2 \notag \\ 
&+ 2 \sum_{p,q=1}^n a_{pq} \, b_{pq} - 2 \sum_{p,q=1}^n e_{pq} \, f_{pq} - 2 \sum_{p,q=1}^n c_{pq} \, c_{qp} - 2 \sum_{p,q=1}^n d_{pq} \, d_{qp}. \notag 
\end{align} 
The assertion follows immediately from (\ref{ingredient.1}) and (\ref{ingredient.2}).

\section{Proof of Theorem \ref{Harnack.inequality}}

We define a Riemannian metric $h$ on $M \times (0,T)$ by 
\[h = \sum_{i,j=1}^n g_{ij} \, dx^i \otimes dx^j + \frac{1}{t^2} \, dt \otimes dt.\]

\begin{lemma}
\label{estimates.for.h}
Suppose that 
\[\sup_{(x,t) \in M \times (0,T)} |\text{\rm Rm}| < \infty.\] 
Then there exists a uniform constant $C$ such that 
\[\Big | \tilde{D}_{\frac{\partial}{\partial t}} h - \tilde{\Delta} h - \frac{1}{t} \, h \Big |_h \leq C\] 
and 
\[|\tilde{D}_v h|_h \leq C \, |v|\] 
for all points $(x,t) \in M \times (0,T)$ and all vectors $v \in T_x M$.
\end{lemma}

\textbf{Proof.} 
By definition of $\tilde{D}$, we have 
\begin{align*} 
&\tilde{D}_{\frac{\partial}{\partial x^i}} dx^j = -\Gamma_{ik}^j \, dx^k + \Big ( \text{\rm Ric}_i^j + \frac{1}{2t} \, \delta_i^j \Big ) \, dt \\ 
&\tilde{D}_{\frac{\partial}{\partial x^i}} dt = 0 \\ 
&\tilde{D}_{\frac{\partial}{\partial t}} dx^j = \Big ( \text{\rm Ric}_i^j + \frac{1}{2t} \, \delta_i^j \Big ) \, dx^i + \frac{1}{2} \, \partial^j \text{\rm scal} \, dt \\ 
&\tilde{D}_{\frac{\partial}{\partial t}} dt = \frac{3}{2t} \, dt. 
\end{align*} 
This implies 
\[\tilde{D}_{\frac{\partial}{\partial x^k}} h = \Big ( \text{\rm Ric}_{ik} + \frac{1}{2t} \, g_{ik} \Big ) \, (dx^i \otimes dt + dt \otimes dx^i).\] 
Moreover, we have 
\begin{align*} 
\tilde{D}_{\frac{\partial}{\partial t}} h 
&= \frac{1}{2} \, \partial_i \text{\rm scal} \, (dx^i \otimes dt + dt \otimes dx^i) \\ 
&+ \frac{1}{t} \, g_{ij} \, dx^i \otimes dx^j + \frac{1}{t^3} \, dt \otimes dt 
\end{align*} 
and 
\begin{align*} 
\tilde{\Delta} h 
&= \frac{1}{2} \, \partial_i \text{\rm scal} \, (dx^i \otimes dt + dt \otimes dx^i) \\ 
&+ 2 \, \Big ( \text{\rm Ric}_k^i + \frac{1}{2t} \, \delta_k^i \Big ) \, \Big ( \text{\rm Ric}_i^k + \frac{1}{2t} \, \delta_i^k \Big ) \, dt \otimes dt. 
\end{align*}
Putting these facts together, we obtain 
\[\tilde{D}_{\frac{\partial}{\partial t}} h - \tilde{\Delta} h - \frac{1}{t} \, h 
= -2 \, \Big ( \text{\rm Ric}_k^i + \frac{1}{2t} \, \delta_k^i \Big ) \, \Big ( \text{\rm Ric}_i^k + \frac{1}{2t} \, \delta_i^k \Big ) \, dt \otimes dt.\] 
Thus, we conclude that 
\[\Big | \tilde{D}_{\frac{\partial}{\partial t}} h - \tilde{\Delta} h - \frac{1}{t} \, h \Big |_h = 2t^2 \, |\text{\rm Ric}|^2 + 2t \, \text{\rm scal} + \frac{n}{2}\] 
and 
\[|\tilde{D}_v h|_h^2 = 2t^2 \, \text{\rm Ric}^2(v,v) + 2t \, \text{\rm Ric}(v,v) + \frac{1}{2} \, g(v,v)\] 
for all points $(x,t) \in M \times (0,T)$ and all vectors $v \in T_x M$. \\

\begin{lemma}
\label{preparation}
Suppose that $(M,g(t)) \times \mathbb{R}^2$ has nonnegative isotropic curvature for all $t \in (0,T)$. Moreover, we assume that 
\[\sup_{(x,t) \in M \times (0,T)} |D^m \text{\rm Rm}| < \infty\] 
for $m = 0,1,2,\hdots$. Then there exists a uniform constant $C$ such that 
\[S + \frac{1}{4} \, C \, t \: h \owedge h \in K\] 
for all points $(x,t) \in M \times (0,T)$. Here, $\owedge$ denotes the Kulkarni-Nomizu product.
\end{lemma}

\textbf{Proof.} 
There exists a uniform constant $C$ such that 
\[\big | S - R_{ijkl} \, dx^i \otimes dx^j \otimes dx^k \otimes dx^l \big |_h \leq C \, t\] 
for all $(x,t) \in M \times (0,T)$. This implies 
\[S - R_{ijkl} \, dx^i \otimes dx^j \otimes dx^k \otimes dx^l + \frac{1}{4} \, C \, t \: h \owedge h \in K\] 
for all $(x,t) \in M \times (0,T)$. Moreover, since $(M,g(t)) \times \mathbb{R}^2$ has nonnegative isotropic curvature, we have 
\[R_{ijkl} \, dx^i \otimes dx^j \otimes dx^k \otimes dx^l \in K\] 
for all points $(x,t) \in M \times (0,T)$. Putting these facts together, the assertion follows. \\

\begin{proposition}
\label{key.step}
Suppose that $(M,g(t)) \times \mathbb{R}^2$ has nonnegative isotropic curvature for all $t \in (0,T)$. Moreover, we assume that 
\[\sup_{(x,t) \in M \times (0,T)} |D^m \text{\rm Rm}| < \infty\] 
for $m = 0,1,2,\hdots$. Then $S_{(x,t)} \in K$ for all $(x,t) \in M \times (0,T)$. 
\end{proposition}

\textbf{Proof.} 
By Lemma 5.1 in \cite{Hamilton4}, we can find a smooth function $\varphi: M \to \mathbb{R}$ with the following properties: 
\begin{itemize}
\item[(i)] $\varphi(x) \to \infty$ as $x \to \infty$ 
\item[(ii)] $\varphi(x) \geq 1$ for all $x \in M$ 
\item[(iii)] $\sup_{(x,t) \in M \times (0,T)} |\nabla \varphi(x)|_{g(t)} < \infty$
\item[(iv)] $\sup_{(x,t) \in M \times (0,T)} |\Delta_{g(t)} \varphi(x)| < \infty$ 
\end{itemize} 
Let $\varepsilon$ be an arbitrary positive real number. We define a $(0,4)$-tensor $\hat{S}$ by 
\[\hat{S} = S + \frac{1}{4} \, \varepsilon \, e^{\lambda t} \, \varphi(x) \, h \owedge h,\] 
where $\lambda$ is a positive constant that will be specified later. Clearly, $\hat{S}$ is an 
algebraic curvature tensor. By Lemma \ref{preparation}, there exists a uniform constant $C_1$ such that 
\[S + \frac{1}{4} \, C_1 \, t \: h \owedge h \in K\] 
for all points $(x,t) \in M \times (0,T)$. Hence, if $\varepsilon \, e^{\lambda t} \, \varphi(x) > C_1 \, t$, 
then $\hat{S}_{(x,t)}$ lies in the interior of the cone $K$.

We claim that $\hat{S}_{(x,t)} \in K$ for all $(x,t) \in M \times (0,T)$. 
Suppose this is false. Then there exists a point $(x_0,t_0) \in M \times (0,T)$ such that 
$\hat{S}_{(x_0,t_0)} \in \partial K$ and $\hat{S}_{(x,t)} \in K$ for all 
$(x,t) \in M \times (0,t_0]$. Since $\hat{S}_{(x_0,t_0)} \in \partial K$, we can find vectors $\tilde{v}_1,\tilde{v}_2,\tilde{v}_3,\tilde{v}_4 \in T_{(x_0,t_0)}(M \times (0,T))$ 
such that 
\[|\tilde{v}_1 \wedge \tilde{v}_3 + \tilde{v}_4 \wedge \tilde{v}_2|_h^2 + |\tilde{v}_1 \wedge \tilde{v}_4 + \tilde{v}_2 \wedge \tilde{v}_3|_h^2 > 0\] 
and 
\begin{align*} 
&\hat{S}(\tilde{v}_1,\tilde{v}_3,\tilde{v}_1,\tilde{v}_3) + \hat{S}(\tilde{v}_1,\tilde{v}_4,\tilde{v}_1,\tilde{v}_4) \\ 
&+ \hat{S}(\tilde{v}_2,\tilde{v}_3,\tilde{v}_2,\tilde{v}_3) + \hat{S}(\tilde{v}_2,\tilde{v}_4,\tilde{v}_2,\tilde{v}_4) - 2 \, \hat{S}(\tilde{v}_1,\tilde{v}_2,\tilde{v}_3,\tilde{v}_4) = 0 
\end{align*} 
at $(x_0,t_0)$. It follows from Proposition \ref{algebraic.fact} that 
\begin{align}
\label{estimate.1}
&\tilde{Q}(\hat{S})(\tilde{v}_1,\tilde{v}_3,\tilde{v}_1,\tilde{v}_3) - \tilde{Q}(\hat{S})(\tilde{v}_1,\tilde{v}_4,\tilde{v}_1,\tilde{v}_4) \notag \\ 
&+ \tilde{Q}(\hat{S})(\tilde{v}_2,\tilde{v}_3,\tilde{v}_2,\tilde{v}_3) - \tilde{Q}(\hat{S})(\tilde{v}_2,\tilde{v}_4,\tilde{v}_2,\tilde{v}_4) + 2 \, \tilde{Q}(\hat{S})(\tilde{v}_1,\tilde{v}_2,\tilde{v}_3,\tilde{v}_4) \geq 0
\end{align} 
at $(x_0,t_0)$. We may extend $\tilde{v}_1,\tilde{v}_2,\tilde{v}_3,\tilde{v}_4$ to vector fields on $M \times (0,T)$ such that 
\begin{align*} 
&\tilde{D}_{\frac{\partial}{\partial x^i}} \tilde{v}_1 = 0 \qquad\qquad \tilde{D}_{\frac{\partial}{\partial t}} \tilde{v}_1 - \tilde{\Delta} \tilde{v}_1 + \frac{1}{2t} \, \tilde{v}_1 = 0 \\
&\tilde{D}_{\frac{\partial}{\partial x^i}} \tilde{v}_2 = 0 \qquad\qquad \tilde{D}_{\frac{\partial}{\partial t}} \tilde{v}_2 - \tilde{\Delta} \tilde{v}_2 + \frac{1}{2t} \, \tilde{v}_2 = 0 \\
&\tilde{D}_{\frac{\partial}{\partial x^i}} \tilde{v}_3 = 0 \qquad\qquad \tilde{D}_{\frac{\partial}{\partial t}} \tilde{v}_3 - \tilde{\Delta} \tilde{v}_3 + \frac{1}{2t} \, \tilde{v}_3 = 0 \\
&\tilde{D}_{\frac{\partial}{\partial x^i}} \tilde{v}_4 = 0 \qquad\qquad \tilde{D}_{\frac{\partial}{\partial t}} \tilde{v}_4 - \tilde{\Delta} \tilde{v}_4 + \frac{1}{2t} \, \tilde{v}_4 = 0
\end{align*} 
at $(x_0,t_0)$. We now define a function $f: M \times (0,T) \to \mathbb{R}$ by 
\begin{align*} 
f &= \hat{S}(\tilde{v}_1,\tilde{v}_3,\tilde{v}_1,\tilde{v}_3) + \hat{S}(\tilde{v}_1,\tilde{v}_4,\tilde{v}_1,\tilde{v}_4) \\ 
&+ \hat{S}(\tilde{v}_2,\tilde{v}_3,\tilde{v}_2,\tilde{v}_3) + \hat{S}(\tilde{v}_2,\tilde{v}_4,\tilde{v}_2,\tilde{v}_4) - 2 \, \hat{S}(\tilde{v}_1,\tilde{v}_2,\tilde{v}_3,\tilde{v}_4). 
\end{align*} 
Clearly, $f(x_0,t_0) = 0$ and $f(x,t) \geq 0$ for all $(x,t) \in M \times (0,t_0]$. This implies 
\[\frac{\partial}{\partial t} f - \Delta f \leq 0\] 
at $(x_0,t_0)$. Hence, if we put 
\[Z = \tilde{D}_{\frac{\partial}{\partial t}} \hat{S} - \tilde{\Delta} \hat{S} - \frac{2}{t} \, \hat{S},\] 
then we obtain 
\begin{align}
\label{estimate.2}
&Z(\tilde{v}_1,\tilde{v}_3,\tilde{v}_1,\tilde{v}_3) + Z(\tilde{v}_1,\tilde{v}_4,\tilde{v}_1,\tilde{v}_4) \notag \\ 
&+ Z(\tilde{v}_2,\tilde{v}_3,\tilde{v}_2,\tilde{v}_3) + Z(\tilde{v}_2,\tilde{v}_4,\tilde{v}_2,\tilde{v}_4) - 2 \, Z(\tilde{v}_1,\tilde{v}_2,\tilde{v}_3,\tilde{v}_4) \\ 
&= \frac{\partial}{\partial t} f - \Delta f \leq 0 \notag 
\end{align} 
at $(x_0,t_0)$. On the other hand, it follows from Proposition \ref{evolution.equation} that
\begin{align*} 
Z &= \tilde{D}_{\frac{\partial}{\partial t}} \hat{S} - \tilde{\Delta} \hat{S} - \frac{2}{t} \, \hat{S} \\ 
&= \tilde{Q}(S) + \frac{1}{4} \, \lambda \, \varepsilon \, e^{\lambda t} \, \varphi(x) \, h \owedge h 
- \frac{1}{4} \, \varepsilon \, e^{\lambda t} \, \Delta \varphi(x) \, h \owedge h \\ 
&- \varepsilon \, e^{\lambda t} \, \sum_{j=1}^n \langle \nabla \varphi(x),e_j \rangle \: h \owedge \tilde{D}_{e_j} h \\ 
&- \frac{1}{2} \, \varepsilon \, e^{\lambda t} \, \varphi(x) \, \sum_{j=1}^n \tilde{D}_{e_j} h \owedge \tilde{D}_{e_j} h \\ 
&+ \frac{1}{2} \, \varepsilon \, e^{\lambda t} \, \varphi(x) \, h \owedge \Big ( \tilde{D}_{\frac{\partial}{\partial t}} h - \tilde{\Delta} h - \frac{1}{t} \, h \Big ) 
\end{align*} 
for all $(x,t) \in M \times (0,T)$. In view of Lemma \ref{estimates.for.h}, there exists a uniform constant $C_2$ such that 
\[\Big | Z - \tilde{Q}(S) - \frac{1}{4} \, \lambda \, \varepsilon \, e^{\lambda t} \, \varphi(x) \, h \owedge h \Big |_h \leq C_2 \, \varepsilon \, e^{\lambda t} \, (\varphi(x) + |\nabla \varphi(x)| + |\Delta \varphi(x)|)\] 
for all $(x,t) \in M \times (0,T)$. Since $\nabla \varphi(x)$ and $\Delta \varphi(x)$ are uniformly bounded, it follows that 
\[\Big | Z - \tilde{Q}(S) - \frac{1}{4} \, \lambda \, \varepsilon \, e^{\lambda t} \, \varphi(x) \, h \owedge h \Big |_h 
\leq C_3 \, \varepsilon \, e^{\lambda t} \, \varphi(x)\] 
for all $(x,t) \in M \times (0,T)$. 

We next observe that $\varepsilon \, e^{\lambda t_0} \, \varphi(x_0) \leq C_1 \, t_0$. (Indeed, if 
$\varepsilon \, e^{\lambda t_0} \, \varphi(x_0) < C_1 \, t_0$, then $\hat{S}_{(x_0,t_0)}$ would lie in the interior of the cone $K$, contrary to our choice of $(x_0,t_0)$.) Hence, there exists a uniform constant $C_4$ such that
\[|S|_h + |\hat{S} - S|_h \leq C_4\]
at $(x_0,t_0)$. This implies
\begin{align*} 
|\tilde{Q}(\hat{S}) - \tilde{Q}(S)|_h 
&\leq C_5 \, \big ( |S|_h \, |\hat{S} - S|_h + |\hat{S} - S|_h^2 \big ) \\ 
&\leq C_5 \, C_4 \, |\hat{S} - S|_h \\ 
&\leq C_6 \, \varepsilon \, e^{\lambda t} \, \varphi(x) 
\end{align*} 
at $(x_0,t_0)$. Putting these facts together, we obtain 
\[\Big | Z - \tilde{Q}(\hat{S}) - \frac{1}{4} \, \lambda \, \varepsilon \, e^{\lambda t} \, \varphi(x) \, h \owedge h \Big |_h \leq C_7 \, \varepsilon \, e^{\lambda t} \, \varphi(x)\] 
at $(x_0,t_0)$. This implies 
\[Z - \tilde{Q}(\hat{S}) - \frac{1}{4} \, (\lambda - C_7) \, \varepsilon \, e^{\lambda t} \, \varphi(x) \, h \owedge h \in K\] 
at $(x_0,t_0)$. Hence, if we choose $\lambda > C_7$, then we have
\begin{align} 
\label{estimate.3}
&Z(\tilde{v}_1,\tilde{v}_3,\tilde{v}_1,\tilde{v}_3) + Z(\tilde{v}_1,\tilde{v}_4,\tilde{v}_1,\tilde{v}_4) \notag \\ 
&+ Z(\tilde{v}_2,\tilde{v}_3,\tilde{v}_2,\tilde{v}_3) + Z(\tilde{v}_2,\tilde{v}_4,\tilde{v}_2,\tilde{v}_4) - 2 \, Z(\tilde{v}_1,\tilde{v}_2,\tilde{v}_3,\tilde{v}_4) \notag \\ 
&- \tilde{Q}(\hat{S})(\tilde{v}_1,\tilde{v}_3,\tilde{v}_1,\tilde{v}_3) - \tilde{Q}(\hat{S})(\tilde{v}_1,\tilde{v}_4,\tilde{v}_1,\tilde{v}_4) \\ 
&- \tilde{Q}(\hat{S})(\tilde{v}_2,\tilde{v}_3,\tilde{v}_2,\tilde{v}_3) - \tilde{Q}(\hat{S})(\tilde{v}_2,\tilde{v}_4,\tilde{v}_2,\tilde{v}_4) + 2 \, \tilde{Q}(\hat{S})(\tilde{v}_1,\tilde{v}_2,\tilde{v}_3,\tilde{v}_4) \notag \\ 
&> 0 \notag 
\end{align} 
at $(x_0,t_0)$. The inequality (\ref{estimate.3}) is inconsistent with (\ref{estimate.1}) and (\ref{estimate.2}).
Consequently, we have $\hat{S}_{(x,t)} \in K$ for all points $(x,t) \in M \times (0,T)$.
Since $\varepsilon > 0$ is arbitrary, it follows that $S_{(x,t)} \in K$ for all points $(x,t) \in M \times (0,T)$. \\

\begin{proposition}
\label{main.result}
Suppose that $(M,g(t)) \times \mathbb{R}^2$ has nonnegative isotropic curvature for all $t \in (0,T)$. 
Moreover, we assume that 
\[\sup_{(x,t) \in M \times (\alpha,T)} \text{\rm scal}(x,t) < \infty\] 
for all $\alpha \in (0,T)$. Then $S_{(x,t)} \in K$ for all $(x,t) \in M \times (0,T)$. 
\end{proposition}

\textbf{Proof.} 
Fix a real number $\alpha \in (0,T)$. By assumption, we have 
\[\sup_{(x,t) \in M \times (\alpha,T)} |\text{\rm Rm}| < \infty.\] 
Using Shi's interior derivative estimates, we obtain 
\[\sup_{(x,t) \in M \times (\alpha,T)} |D^m \text{\rm Rm}| < \infty\] 
for $m = 1,2, \hdots$ (see e.g. \cite{Hamilton-survey}, Theorem 13.1). Hence, we can apply Proposition \ref{key.step} to the metrics $g(t + \alpha)$, $t \in (0,T - \alpha)$. Taking the limit as $\alpha \to 0$, the assertion follows. \\

Theorem \ref{Harnack.inequality} is an immediate consequence of Proposition \ref{main.result}. To see this, we consider a point $(x,t) \in M \times (0,T)$ and vectors $v,w \in T_x M$. By Proposition \ref{main.result}, we have 
\begin{align*} 
&S(\tilde{v}_1,\tilde{v}_3,\tilde{v}_1,\tilde{v}_3) + S(\tilde{v}_1,\tilde{v}_4,\tilde{v}_1,\tilde{v}_4) \\ 
&+ S(\tilde{v}_2,\tilde{v}_3,\tilde{v}_2,\tilde{v}_3) + S(\tilde{v}_2,\tilde{v}_4,\tilde{v}_2,\tilde{v}_4) - 2 \, S(\tilde{v}_1,\tilde{v}_2,\tilde{v}_3,\tilde{v}_4) \geq 0 
\end{align*}
for all vectors $\tilde{v}_1,\tilde{v}_2,\tilde{v}_3,\tilde{v}_4 \in T_{(x,t)} (M \times (0,T))$. Hence, if we put 
\[\tilde{v}_1 = \frac{\partial}{\partial t} + v, \qquad \tilde{v}_2 = 0, \qquad \tilde{v}_3 = w, \qquad \tilde{v}_4 = 0,\] 
then we obtain 
\[M(w,w) + 2 \, P(v,w,w) + R(v,w,v,w) \geq 0.\] 
This completes the proof of Theorem \ref{Harnack.inequality}. In order to prove Corollary \ref{trace.Harnack.inequality}, we take the trace over $w$. This yields 
\[\Delta \text{\rm scal} + 2 \, |\text{\rm Ric}|^2 + \frac{1}{t} \, 
\text{\rm scal} + 2 \, \partial_i \text{\rm scal} \, v^i + 2 \, \text{\rm Ric}(v,v) \geq 0.\] 
Hence, Corollary \ref{trace.Harnack.inequality} follows from the identity 
$\frac{\partial}{\partial t} \text{\rm scal} = \Delta \text{\rm scal} + 2 \, |\text{\rm Ric}|^2$. \\

\section{The equality case in the Harnack inequality} 

In this section, we analyze the equality case in the Harnack inequality. Let $(M,g(t))$, $t \in (0,T)$, be a family of complete Riemannian manifolds evolving under Ricci flow. As above, we assume that $(M,g(t)) \times \mathbb{R}^2$ has nonnegative isotropic curvature for all $t \in (0,T)$. Moreover, we require that 
\[\sup_{(x,t) \in M \times (\alpha,T)} \text{\rm scal}(x,t) < \infty\] 
for all $\alpha \in (0,T)$.

Let $E$ be the tangent bundle of $M \times (0,T)$. We denote by $P$ the total space of the vector 
bundle $E \oplus E \oplus E \oplus E$. The connection $\tilde{D}$ defines a horizontal distribution on $P$. Hence, the tangent bundle of $P$ splits as a direct sum $TP = \mathbb{H} \oplus \mathbb{V}$, where $\mathbb{H}$ 
and $\mathbb{V}$ denote the horizontal and vertical distributions, respectively. 

Let $\pi$ be the projection from $P$ to $M \times (0,T)$. For each $t \in (0,T)$, we 
denote by $P_t = \pi^{-1}(M \times \{t\})$ the time $t$ slice of $P$. We define a function 
$u: P \to \mathbb{R}$ by 
\begin{align*} 
u: (\tilde{v}_1,\tilde{v}_2,\tilde{v}_3,\tilde{v}_4) 
\mapsto \; &S(\tilde{v}_1,\tilde{v}_3,\tilde{v}_1,\tilde{v}_3) + S(\tilde{v}_1,\tilde{v}_4,\tilde{v}_1,\tilde{v}_4) \\ 
&+ S(\tilde{v}_2,\tilde{v}_3,\tilde{v}_2,\tilde{v}_3) + S(\tilde{v}_2,\tilde{v}_4,\tilde{v}_2,\tilde{v}_4) \\ &- 2 \, S(\tilde{v}_1,\tilde{v}_2,\tilde{v}_3,\tilde{v}_4). 
\end{align*} 
By Proposition \ref{main.result}, $u$ is a nonnegative function on $P$. Let $F = \{u = 0\}$ be the zero set 
of the function $u$. We claim that $F$ is invariant under parallel transport:

\begin{proposition} 
\label{F.is.invariant.under.parallel.transport}
Fix a real number $t_0 \in (0,T)$, and let $\tilde{\gamma}: [0,1] \to P_{t_0}$ be a smooth horizontal 
curve such that $\tilde{\gamma}(0) \in F$. Then $\tilde{\gamma}(s) \in F$ for all $s \in [0,1]$.
\end{proposition}

\textbf{Proof.} 
Without loss of generality, we may assume that the projected path $\pi \circ \tilde{\gamma}: [0,1] \to M \times \{t_0\}$ is contained in a single coordinate chart. 
Let $\Omega \subset M \times (0,T)$ be a coordinate chart such that $\pi(\tilde{\gamma}(s)) \in \Omega$ 
for all $s \in [0,1]$. We can find smooth vector fields $X_1, \hdots, X_n$ on $\Omega$ such that 
\[\sum_{k=1}^n X_k \otimes X_k = \sum_{i,j=1}^n g^{ij} \, \frac{\partial}{\partial x^i} \otimes \frac{\partial}{\partial x^j}.\] 
Moreover, we define a vector field $Y$ on $\Omega$ by 
\[Y = \frac{\partial}{\partial t} + \sum_{k=1}^n \tilde{D}_{X_k} X_k.\] 
Let $\tilde{X}_1,\hdots,\tilde{X}_n,\tilde{Y}$ be the horizontal lifts of $X_1,\hdots,X_n,Y$. 
At each point $(\tilde{v}_1,\tilde{v}_2,\tilde{v}_3,\tilde{v}_4) \in \pi^{-1}(\Omega)$, we have 
\begin{align*} 
&\tilde{Y}(u) - \sum_{k=1}^n \tilde{X}_k(\tilde{X}_k(u)) \\ 
&= \Big ( \tilde{D}_{\frac{\partial}{\partial t}} S - \tilde{\Delta} S \Big ) (\tilde{v}_1,\tilde{v}_3,\tilde{v}_1,\tilde{v}_3) 
+ \Big ( \tilde{D}_{\frac{\partial}{\partial t}} S - \tilde{\Delta} S \Big ) (\tilde{v}_1,\tilde{v}_4,\tilde{v}_1,\tilde{v}_4) \\ 
&+ \Big ( \tilde{D}_{\frac{\partial}{\partial t}} S - \tilde{\Delta} S \Big ) (\tilde{v}_2,\tilde{v}_3,\tilde{v}_2,\tilde{v}_3) 
+ \Big ( \tilde{D}_{\frac{\partial}{\partial t}} S - \tilde{\Delta} S \Big ) (\tilde{v}_2,\tilde{v}_4,\tilde{v}_2,\tilde{v}_4) \\ 
&- 2 \, \Big ( \tilde{D}_{\frac{\partial}{\partial t}} S - \tilde{\Delta} S \Big ) (\tilde{v}_1,\tilde{v}_2,\tilde{v}_3,\tilde{v}_4). 
\end{align*}
Using Proposition \ref{evolution.equation}, we obtain 
\begin{align*} 
&\tilde{Y}(u) - \sum_{k=1}^n \tilde{X}_k(\tilde{X}_k(u)) - \frac{2}{t} \, u \\ 
&= \tilde{Q}(S)(\tilde{v}_1,\tilde{v}_3,\tilde{v}_1,\tilde{v}_3) 
+ \tilde{Q}(S)(\tilde{v}_1,\tilde{v}_4,\tilde{v}_1,\tilde{v}_4) \\ 
&+ \tilde{Q}(S)(\tilde{v}_2,\tilde{v}_3,\tilde{v}_2,\tilde{v}_3) 
+ \tilde{Q}(S)(\tilde{v}_2,\tilde{v}_4,\tilde{v}_2,\tilde{v}_4) \\ 
&- 2 \, \tilde{Q}(S)(\tilde{v}_1,\tilde{v}_2,\tilde{v}_3,\tilde{v}_4) 
\end{align*} 
for all points $(\tilde{v}_1,\tilde{v}_2,\tilde{v}_3,\tilde{v}_4) \in \pi^{-1}(\Omega)$. 
Moreover, it follows from the calculations in Section 3 that 
\begin{align*}
&\tilde{Q}(S)(\tilde{v}_1,\tilde{v}_3,\tilde{v}_1,\tilde{v}_3) + \tilde{Q}(S)(\tilde{v}_1,\tilde{v}_4,\tilde{v}_1,\tilde{v}_4) \\ 
&+ \tilde{Q}(S)(\tilde{v}_2,\tilde{v}_3,\tilde{v}_2,\tilde{v}_3) 
+ \tilde{Q}(S)(\tilde{v}_2,\tilde{v}_4,\tilde{v}_2,\tilde{v}_4) - 2 \, \tilde{Q}(S)(\tilde{v}_1,\tilde{v}_2,\tilde{v}_3,\tilde{v}_4) \\ 
&\geq C \, \inf_{\xi \in \mathbb{V}, \, |\xi| \leq 1} (D^2 u)(\xi,\xi) 
\end{align*} 
for all points $(\tilde{v}_1,\tilde{v}_2,\tilde{v}_3,\tilde{v}_4) \in \pi^{-1}(\Omega)$. 
Here, $D^2 u$ denotes the Hessian of $u$ in vertical direction. Putting these facts together, we obtain 
\[\tilde{Y}(u) - \sum_{k=1}^n \tilde{X}_k(\tilde{X}_k(u)) - \frac{2}{t} \, u \geq C \, \inf_{\xi \in \mathbb{V}, \, |\xi| \leq 1} (D^2 u)(\xi,\xi)\] 
on $\pi^{-1}(\Omega)$. Hence, the assertion follows from J.M.~Bony's version of the strong maximum principle (see \cite{Bony} or \cite{Brendle-Schoen2}, Proposition 4). \\

For each point $(x,t) \in M \times (0,T)$, we denote by $\mathcal{N}_{(x,t)}$ the set of all vectors of the form $\tilde{v} 
= \frac{\partial}{\partial t} + v \in T_{(x,t)}(M \times (0,T))$, where $v \in T_x M$ satisfies 
\[\frac{\partial}{\partial t} \text{\rm scal} + \frac{1}{t} \, \text{\rm scal} + 2 \, \partial_i \text{\rm scal} \, v^i + 2 \, \text{\rm Ric}(v,v) = 0.\] 
In view of Theorem \ref{Harnack.inequality}, we can characterize the set $\mathcal{N}_{(x,t)}$ as follows: 
\begin{align*}
&\frac{\partial}{\partial t} + v \in \mathcal{N}_{(x,t)} \\ 
&\Longleftrightarrow \quad \frac{\partial}{\partial t} \text{\rm scal} + \frac{1}{t} \, \text{\rm scal} + 2 \, \partial_i \text{\rm scal} \, v^i + 2 \, \text{\rm Ric}(v,v) = 0 \\ 
&\Longleftrightarrow \quad M(w,w) + 2 \, P(v,w,w) + R(v,w,v,w) = 0 \quad \text{\rm for all $w \in T_x M$} \\ 
&\Longleftrightarrow \quad \Big ( \frac{\partial}{\partial t} + v,0,w,0 \Big ) \in F \quad \text{\rm for all $w \in T_x M$} 
\end{align*}
By Proposition \ref{F.is.invariant.under.parallel.transport}, the set $F$ is invariant under parallel transport. Therefore, we can draw the following conclusion:

\begin{corollary}
\label{equality.case.in.the.trace.Harnack.inequality}
Fix a smooth path $\gamma: [0,1] \to M \times \{t_0\}$. We denote by $\tilde{P}_\gamma: T_{\gamma(0)}(M \times (0,T)) \to T_{\gamma(1)}(M \times (0,T))$ the parallel transport along $\gamma$ with respect to the connection $\tilde{D}$. If 
$\tilde{v} \in \mathcal{N}_{\gamma(0)}$, then $\tilde{P}_\gamma \tilde{v} \in \mathcal{N}_{\gamma(1)}$.
\end{corollary}

\begin{proposition}
\label{expanding.soliton}
Let $(M,g(t))$, $t \in (0,T)$, be a family of complete Riemannian manifolds evolving under Ricci flow. For each $t \in (0,T)$, we assume that $(M,g(t)) \times \mathbb{R}^2$ has nonnegative isotropic curvature and $(M,g(t))$ has positive Ricci curvature. Moreover, suppose that there exists a point $(x_0,t_0) \in M \times (0,T)$ such that 
\[t_0 \cdot \text{\rm scal}(x_0,t_0) = \sup_{(x,t) \in M \times (0,T)} t \cdot \text{\rm scal}(x,t).\] 
Then there exists a smooth vector field $V = V^j \, \frac{\partial}{\partial x^j}$ such that 
\[D_{\frac{\partial}{\partial x^i}} V = \text{\rm Ric}_i^j \, \frac{\partial}{\partial x^j} + \frac{1}{2t} \, \frac{\partial}{\partial x^i}\] 
for all $(x,t) \in M \times \{t_0\}$. In particular, $(M,g(t_0))$ is an expanding Ricci soliton.
\end{proposition}

\textbf{Proof.} 
Since $(M,g(t))$ has positive Ricci curvature, there exists a unique vector field $V = V^j \, \frac{\partial}{\partial x^j}$ such that $\partial_i \text{\rm scal} + 2 \, \text{\rm Ric}_{ij} \, V^j = 0$. We claim that 
\begin{equation} 
\label{inclusion}
\mathcal{N}_{(x,t)} \subset \Big \{ \frac{\partial}{\partial t} + V_{(x,t)} \Big \} 
\end{equation}
for all points $(x,t) \in M \times (0,T)$. In order to prove this, we consider an arbitrary vector 
$\tilde{v} \in \mathcal{N}_{(x,t)}$. The vector $\tilde{v}$ can be written in the form $\tilde{v} = \frac{\partial}{\partial t} + v$, where $v \in T_x M$ satisfies 
\[\frac{\partial}{\partial t} \text{\rm scal} + \frac{1}{t} \, \text{\rm scal} + 2 \, \partial_i \text{\rm scal} \, v^i + 2 \, \text{\rm Ric}(v,v) = 0.\] 
Using Corollary \ref{trace.Harnack.inequality}, we conclude that $\partial_i \text{\rm scal} + 2 \, \text{\rm Ric}_{ij} \, v^j = 0$. Since $(M,g(t))$ has positive Ricci curvature, it follows that $v = V_{(x,t)}$. This completes the proof of (\ref{inclusion}). In particular, the set $\mathcal{N}_{(x,t)}$ contains at most one element. 

By assumption, the function $t \cdot \text{\rm scal}(x,t)$ attains its global maximum at $(x_0,t_0)$. This implies 
\[\frac{\partial}{\partial t} \text{\rm scal} + \frac{1}{t} \, \text{\rm scal} = 0\] 
at $(x_0,t_0)$. Consequently, the set $\mathcal{N}_{(x_0,t_0)}$ is non-empty. Hence, it follows from Corollary \ref{equality.case.in.the.trace.Harnack.inequality} that the set $\mathcal{N}_{(x,t)}$ is non-empty for all points $(x,t) \in M \times \{t_0\}$. Using (\ref{inclusion}), we obtain 
\begin{equation} 
\label{identity}
\mathcal{N}_{(x,t)} = \Big \{ \frac{\partial}{\partial t} + V_{(x,t)} \Big \}
\end{equation} 
for all points $(x,t) \in M \times \{t_0\}$. Hence, by Corollary \ref{equality.case.in.the.trace.Harnack.inequality}, we have 
\[\tilde{P}_\gamma \Big ( \frac{\partial}{\partial t} + V_{\gamma(0)} \Big ) = \frac{\partial}{\partial t} + V_{\gamma(1)}\] 
for every smooth path $\gamma: [0,1] \to M \times \{t_0\}$. Thus, we conclude that 
\[\tilde{D}_{\frac{\partial}{\partial x^i}} \Big ( \frac{\partial}{\partial t} + V \Big ) = 0\] 
for all points $(x,t) \in M \times \{t_0\}$. From this, the assertion follows. \\

\section{Ancient solutions to the Ricci flow}

In this final section, we consider ancient solutions to the Ricci flow. In this case, we are able to remove the $1/t$ terms in the Harnack inequality: 

\begin{proposition} 
Let $(M,g(t))$, $t \in (-\infty,T)$, be a family of complete Riemannian manifolds evolving under Ricci flow. We assume that $(M,g(t)) \times \mathbb{R}^2$ has nonnegative isotropic curvature for all $t \in (-\infty,T)$. Moreover, we assume that \[\sup_{(x,t) \in M \times (\alpha,T)} \text{\rm scal}(x,t) < \infty\] 
for all $\alpha \in (-\infty,T)$. Then we have 
\[\frac{\partial}{\partial t} \text{\rm scal} + 2 \, \partial_i \text{\rm scal} \, v^i + 2 \, \text{\rm Ric}(v,v) \geq 0\] 
for all points $(x,t) \in M \times (-\infty,T)$ and all vectors $v \in T_x M$. 
\end{proposition}

\textbf{Proof.} 
We employ an argument due to R.~Hamilton \cite{Hamilton5}. To that end, we fix a real number $\alpha \in (-\infty,T)$, and apply Corollary \ref{trace.Harnack.inequality} to the metrics $g(t + \alpha)$, $t \in (0,T - \alpha)$. This implies 
\[\frac{\partial}{\partial t} \text{\rm scal} + \frac{1}{t-\alpha} \, \text{\rm scal} + 2 \, \partial_i \text{\rm scal} \, v^i + 2 \, \text{\rm Ric}(v,v) \geq 0\] 
for all points $(x,t) \in M \times (\alpha,T)$ and all $v \in T_x M$. Taking the limit as $\alpha \to -\infty$, the assertion follows. \\

Our last result generalizes Theorem 1.1 in \cite{Hamilton5}: 

\begin{proposition}
\label{steady.soliton}
Let $(M,g(t))$, $t \in (-\infty,T)$, be a family of complete Riemannian manifolds evolving under Ricci flow. For each $t \in (-\infty,T)$, we assume that $(M,g(t)) \times \mathbb{R}^2$ has nonnegative isotropic curvature and $(M,g(t))$ has positive Ricci curvature. Moreover, suppose that there exists a point $(x_0,t_0) \in M \times (-\infty,T)$ such that 
\[\text{\rm scal}(x_0,t_0) = \sup_{(x,t) \in M \times (-\infty,T)} \text{\rm scal}(x,t).\] 
Then there exists a smooth vector field $V = V^j \, \frac{\partial}{\partial x^j}$ such that 
\[D_{\frac{\partial}{\partial x^i}} V = \text{\rm Ric}_i^j \, \frac{\partial}{\partial x^j}\] 
for all $(x,t) \in M \times \{t_0\}$. In particular, $(M,g(t_0))$ is a steady Ricci soliton.
\end{proposition}

The proof of Proposition \ref{steady.soliton} is analogous to the proof of 
Proposition \ref{expanding.soliton} above. The details are left to the reader.


\begin{thebibliography}{99}
\bibitem{Bony} 
J.M.~Bony, \textit{Principe du maximum, in\'egalit\'e de Harnack et unicit\'e du probl\`eme de Cauchy pour les op\'erateurs elliptiques d\'eg\'en\'er\'es,} Ann. Inst. Fourier (Grenoble), 19, 277--304 (1969) 

\bibitem{Brendle-Schoen1} 
S.~Brendle and R.~Schoen, \textit{Manifolds with $1/4$-pinched curvature are space forms,} J. Amer. Math. Soc. (to appear)

\bibitem{Brendle-Schoen2}
S.~Brendle and R.~Schoen, \textit{Classification of manifolds with weakly $1/4$-pinched curvatures,} Acta Math. 200, 1--13 (2008)

\bibitem{Cao}
H.D.~Cao, \textit{On Harnack's inequalities for the K\"ahler-Ricci flow,} Invent. Math. 109, 247--263 (1992)

\bibitem{Chow-Chu} 
B.~Chow and S.C.~Chu, \textit{A geometric interpretation of Hamilton's Harnack inequality for the Ricci flow,} Math. Res. Lett. 2, 701--718 (1995)

\bibitem{Chow-Knopf}
B.~Chow and D.~Knopf, \textit{New Li-Yau-Hamilton inequalities for the Ricci flow via the space-time approach,} J. Diff. Geom. 60, 1--54 (2002)

\bibitem{Hamilton1}
R.~Hamilton, \textit{Three-manifolds with positive Ricci curvature,} J. Diff. Geom. 17, 255--306 (1982)

\bibitem{Hamilton2}
R.~Hamilton, \textit{Four-manifolds with positive curvature operator,} J. Diff. Geom. 24, 153--179 (1986)

\bibitem{Hamilton3}
R.~Hamilton, \textit{The Ricci flow on surfaces,} Contemp. Math. 71, 237--262 (1988)

\bibitem{Hamilton4} 
R.~Hamilton, \textit{The Harnack estimate for the Ricci flow,} J. Diff. Geom. 37, 225--243 (1993)

\bibitem{Hamilton5}
R.~Hamilton, \textit{Eternal solutions to the Ricci flow,} J. Diff. Geom. 38, 1--11 (1993)

\bibitem{Hamilton-survey}
R.~Hamilton, \textit{The formation of singularities in the Ricci flow,} Surveys in Differential Geometry, Vol. II, Intl. Press, Cambridge MA, 1995, 7--136
\end{thebibliography}
\end{document}